\documentclass[12pt]{article}
\usepackage{amsmath,amsthm,amssymb,amsfonts}

\newfont{\lie}{eufm10 at 12pt}
\newfont{\liepequenos}{eufm10 at 10pt}

\newfont{\corpos}{msbm10 at 12pt}
\newfont{\corpospequenos}{msbm10 at 10pt}

\newcommand{\C}{\mbox{\corpos \symbol{67}}}          
\newcommand{\Hamil}{\mbox{\corpos \symbol{72}}}      
\newcommand{\Pro}{\mbox{\corpos \symbol{80}}}        
\newcommand{\R}{\mbox{\corpos \symbol{82}}}          

\newcommand{\g}{\mbox{\lie g}}       %
\newcommand{\m}{\mbox{\lie m}}       %

\newcommand{\s}{\mbox{\lie s}}

\newcommand{\uni}{\mbox{\lie u}}     %

\newcommand{\sol}{\mbox{\lie so}}

\newcommand{\XIS}{\mbox{\lie X}}   

\newcommand{\rr}{\rightarrow}
\newcommand{\lrr}{\longrightarrow}

\newcommand{\rdoisn}{\R^{2n}}                    %
\newcommand{\jnab}{{\cal J}^{\nabla}}             
\newcommand{\hnab}{{\cal H}^{\nabla}}             
\newcommand{\zo}{{\cal Z}}           

\newcommand{\nab}[2]{{\nabla_{#1}{#2}}}
\newcommand{\Tr}[1]{{\mathrm{Tr}}\,{#1}}

\newcommand{\End}[1]{{\mathrm{End}}\,{#1}}

\newcommand{\Id}{{\mathrm{Id}}}
\newcommand{\dx}{{\mathrm{d}}}
\newcommand{\inv}[1]{{#1}^{-1}}
\newcommand{\papa}[2]{\frac{\partial#1}{\partial#2}}

\newcommand{\db}{\overline{\partial}}
\newcommand{\cinf}[1]{{\mathrm{C}}^\infty_{#1}}
\newcommand{\prl}{{\mathrm{Re}}\,}
\newcommand{\pig}{{\mathrm{Im}}\,}

\newcommand{\vol}{{\mathrm{Vol}}\,}

\newtheorem{teo}{Theorem}[section]
\newtheorem{lema}{Lemma}[section]
\newtheorem{coro}{Corollary}[section]
\newtheorem{prop}{Proposition}[section]

\numberwithin{equation}{section}

\def\cyclic{\mathop{\kern0.9ex{{+}
\kern-2.2ex\raise-.28ex\hbox{\Large\hbox{$\circlearrowright$}}}}\limits}

\begin{document}

\title{\textit{On the twistor space of pseudo-spheres}}
\author{R. Albuquerque$^{1*}$ and  Isabel M.C.\ Salavessa$^{2*}$
\\[2mm]}
\protect\footnotetext{
$1*~~$Partially supported by Funda\c{c}\~{a}o Ci\^{e}ncia e Tecnologia 
through POCI/MAT/60671/2004 and the Plurianual of CIMA.\\
$2*~~$Partially supported by Funda\c{c}\~{a}o Ci\^{e}ncia e Tecnologia through 
POCI/MAT/60671/2004 and the Plurianual of CFIF.}
\date{\today}
\maketitle ~~~\\[-5mm]
{\footnotesize $^1$ Centro de Investiga\c c\~ao em Matem\'{a}tica e Aplica\c c\~oes,
Universidade de \'Evora,\\[-1mm]
7000 \'EVORA, Portugal;~~
e-mail: rpa@uevora.pt}\\
{\footnotesize $^2$ Centro de F\'{\i}sica das Interac\c{c}\~{o}es Fundamentais, Instituto Superior
T\'{e}cnico,
Edif\'{\i}cio Ci\^{e}ncia,\\[-1mm] Piso 3,
1049-001 LISBOA, Portugal;~~
e-mail: isabel@cartan.ist.utl.pt}\\[5mm]
{\small {\bf Abstract:} 
We give a new proof that the sphere $S^6$ does not admit an integrable orthogonal complex structure, as in \cite{Lebr}, following the methods from twistor theory. 

We present the twistor space $\zo^{p,q}$ of a pseudo-sphere $S^{2n}_{2q}=SO_{2p+1,2q}/SO_{2p,2q}$ as a pseudo-K\"ahler symmetric space. We then consider orthogonal complex structures on the pseudo-sphere, only to prove such a structure cannot exist.}

\markright{\sl\hfill  Albuquerque--Salavessa \hfill}

\section{Introduction}

This article raises some questions around the problem solved by C.~Lebrun in \cite{Lebr} about the non existence of orthogonal complex structures on the sphere $S^6$. That clever proof recurs to a particular fibre bundle, the open subspace of the Grassmannian $Gr_3(\C^7)$ consisting of 3-planes $P$ for which $P\cap \overline{P}=\{0\}$. This is a space which, we know today, agrees with the general twistor bundle of the 6-sphere.

The reader may notice throughout the text that we somehow reproduce the same arguments from the refered article, but our goal is to present them as a consequence of the theory of twistor spaces. Moreover the final argument is purely geometric rather than topological. We also extend some known results from the Riemannian to the semi-Riemanniann context, for which it is essential to consider all what was explained in \cite{Obri}. In recalling the theory from this reference we are led to some new insights relating affine transformations and the twistor pseudo-holomorphic structure.

In the last section we revise and compute a few metrics on the twistor space of a pseudo-sphere. We start by proving the spheres cannot be pseudo-K\"ahler. Then putting together the pseudo-K\"ahlerian structure of the twistor space and its intrinsic geometry induced by the linear connection, we are able to find interesting formulae dealing with its curvature and a 2-form $\omega$ on the base manifold. This is actually true for all symplectic twistor spaces. 

The analysis of the exterior derivative of the K\"ahler form from two different paths leads to the conclusion that it must vanish. Although an unexpected proof, by a difference in scalars, it may explain why it was not a trivial problem.

\section{Twistor spaces}

Let $(M,\nabla)$ be a $2n$-dimensional manifold endowed with a linear connection. We briefly recall along the text the theory of twistor spaces described in \cite{Obri,Rawnsley}. For a fast exposition and new proofs we avoid mentioning the principal bundle of frames of $M$.

\subsection{The general theory}

Consider the general twistor space of $M$, ie. the bundle 
\begin{equation}
  {\cal J}(M)=\bigl\{j\in \End{T_xM}|\ x\in M,\ j^2=-1\bigr\}\stackrel{\pi}\lrr M
\end{equation}  
with standard fibre $GL_{2n}(\R)/GL_n(\C)$ which consists of the \textit{complex} symmetric space of linear complex structures on $\rdoisn$. More accurately the bundle is called a \textit{twistor} when it is seen with a certain almost complex structure $\jnab$ induced by $\nabla$. First we have an exact sequence of vector bundles (all over the same base space)
\begin{equation}\label{seqexacta}
0\lrr {\cal V}\lrr T{\cal J}(M)\lrr E=\pi^*TM\lrr 0,
\end{equation}
where $\cal V=\ker\dx\pi$. Then we use the connection to find a splitting $T{\cal J}(M)={\cal V}\oplus\hnab$ into vertical and horizontal tangent vectors and define, up to canonical isomorphism $\dx\pi:\hnab\rr E$,
\begin{equation}\label{estrcomplexatwistor}
 \jnab_j(X)=jX,\ \ \mbox{for $X$ horizontal},\quad\quad \jnab_j(A)=jA,\ \ \mbox{for $A$ vertical}.
\end{equation}
The meaning of ``$jA$" on the vertical side is explained as follows. The general twistor's fibre $\inv{\pi}(x)$ consists of elements $j$ of the form $gJ_0\inv{g}$ where $g$ varies in $GL(T_xM)$ and $J_0$ is a fixed element. So it agrees with the complex symmetric space $GL(T_xM)/GL(T_xM,J_0)$. It is not hard to see that
\begin{equation}
 T_j(\inv{\pi}(x))={\cal V}_j=\bigl\{ A\in\End{E_j}|\ Aj=-jA\bigr\}
\end{equation}
and that this space is closed under left multiplication by $j$. This is the symmetric space complex structure of the standard fibre, which we copy to each fibre of the twistor bundle.

If we define a tautological section $\Phi\in\Gamma({\cal J}(M),\End{E})$ by $\Phi_j=j$, then it  varies along the vertical directions only. More precisely:
\begin{prop}[\cite{Obri}]\label{horizontais}
$\hnab=\{X\in T{\cal J}(M)|\ (\pi^*\nabla)_{X}{\Phi}=0\}$. The vertical part of $X\in T{\cal J}(M)$ is $X'=\frac{1}{2}\Phi(\pi^*\nabla_X\Phi)$.
\end{prop}
To see this we may argue with a section $j:U\mapsto {\cal J}(M)$ on a neighborhood $U$ of a point $x_0$. It is well understood that $\dx j_{x_0}(X)$ lies in the horizontal distribution induced by a connection on a fibre bundle if, and only if, $\nab{X_{x_0}}{j}=0$. But immediately we also deduce $(\pi^*\nabla)_{j_*X}{\Phi}=j^*(\pi^*\nabla)_{X}{j^*\Phi}=\nab{X}{j}$. Here is a complete proof of the proposition. Take normal coordinates $x^i$ for $\nabla$ in $M$ around a point $x_0$, so that, if $\nabla=\dx+{\cal A}$, then ${\cal A}_{x_0}=0$. Take coordinates $z^\alpha$ for the fibre of ${\cal J}(M)$ ($\alpha=1,\ldots,n^2-n$). Then at the point $j=(x_0,[z^\alpha])$ the section $\Phi$ corresponds to $[z^\alpha]$, so $\pi^*\nabla_{\partial_i}\Phi=\pi^*(\dx+{\cal A}_{x_0})_{\partial_i}[z^\alpha]=\papa{[z^\alpha]}{x^i}=0$ and 
$(\pi^*\nabla_{\partial_\alpha}\Phi)\partial_i=(\pi^*\dx_{\partial_\alpha}[z^{\beta}])\partial_i =\partial_{\alpha}[z^{\beta}]\partial_i-[z^\beta]\partial_\alpha\partial_i=[\partial_\alpha,\Phi]\partial_i$. Hence, for $A\in \cal V$, we found $\pi^*\nab{A}{\Phi}=[A,\Phi]$.

Now we recall the integrability equations of $\jnab$, the proof being postponed to section 3.3. Let $j^+,j^-$ denote respectively the projections 
\[ \frac{1}{2}(1-ij),\quad\frac{1}{2}(1+ij)  \]
to the $+i$ and $-i$ eigenspaces of $j$.
\begin{teo}[\cite{Obri}]\label{teo2.1}
The twistor space almost complex structure is integrable if and only if the torsion $T$ and the curvature $R$ of $\nabla$ satisfy
\begin{equation}\label{int}
   j^+T(j^-X,j^-Y)=0,\qquad j^+R(j^-X,j^-Y)j^-=0,   
\end{equation}
for all $X,Y\in TM,\ j\in{\cal J}(M)$.
\end{teo}

\subsubsection{The Riemannian twistor space}

When the structure group of $M$ is reducible and $M$ admits a connection compatible with such  reduction, we can further reduce the twistor space. Here are some celebrated examples: for oriented Riemannian manifolds and metric connections the appropriate twistor is the one with fibre $SO_{2n}/U_n$ (cf. \cite{Atiyah,Hit,Obri,Sal} between many others), for almost hermitian manifolds with a hermitian connection one restricts to $U_{p+q}/U_p\times U_q$ (cf. \cite{Bu,Obri}) and for symplectic manifolds endowed with symplectic connections we consider $Sp_{n}(\R)/U_n$ (cf. \cite{AR,Vais}). But some other twistor spaces have been studied, both of the compact and non-compact type. Namely for the quaternionic structure $I,J,K$ in dimension $4n$ one considers the sphere bundle $\{xI+yJ+zK|\ x^2+y^2+z^2=1\}$. As examples of the non-compact type we mention the hyperbolic twistor space, induced by paraquaternionic structures (cf. \cite{Blair}), and the complex structures compatible with a 2-form or $Sp_{p+q}(\R)/U_{p,q}$ case (cf. \cite{Albu,AR}).

Notice all the previous symmetric spaces are complex symmetric subspaces of the whole space of linear complex structures on $\rdoisn$. This follows trivially from the theory in \cite{KobNomi} (as we shall see in a specific case). Hence the integrability equations of all respective twistor spaces are the same as those for the one with general fibre, cf. theorem \ref{teo2.1}.

In case $(M,g)$ is an oriented Riemannian manifold and we consider the first of the previous examples \begin{equation}\label{etpr}
{\cal J}_+(M,g)=\{j\in{\cal J}(M)|\ j^*g=g\ \mbox{and $j$ induces the same orientation}\}
\end{equation}
with the Levi-Civita connection, then it is a well known result in dimension 4 that $\jnab$ is integrable if, and only if, $M$ is self-dual (cf. \cite{Atiyah}). For higher dimensions it was proved in \cite{Obri}, using representation theory, that the integrability equation being satisfied is equivalent to conformal flatness, ie. the vanishing of the Weyl part of the curvature --- which no longer brakes into two irreducibles as it does in 4 dimensions.

We recall the main lines of the proof, which comes from analysis of equation (\ref{int}). Since for all $j$ we have $j^\pm=k(1\pm iJ_0)\inv{k}=kJ_0^\pm\inv{k}$, the curvature condition can be put as $J_0^+\inv{k}R(kJ_0^-X,kJ_0^-Y)kJ_0^-=0$, $\,\forall k\in SO(T_xM),\ X,Y\in T_xM$. Noticing the adjoint action, the condition is saying $R$ takes values in the largest invariant subspace of curvature type tensors which satisfy $J_0^+R(J_0^-X,J_0^-Y)J_0^-=0$. But we may view $J_0$ as an element of the Lie algebra acting by
\[  (J_0\cdot R) (X,Y)=J_0R(X,Y)-R(J_0X,Y)-R(X,J_0Y)-R(X,Y)J_0,   \]
$\forall X,Y\in T_xM$. Since $J_0$ has eigenvalues $\pm i$ on $T_xM$, it can only have $0,\pm2i,\pm4i$ eigenvalues on curvature tensors (a simple computation). The $4i$ eigenspace is easily seen to consist of tensors of the form $J_0^+R(J_0^-X,J_0^-Y)J_0^-$, so, again, the condition is saying $R$ takes values in the largest invariant subspace in which $J_0$ has no $4i$ eigenvalue. By conjugation and since the tensor $R$ is real, we cannot have the $-4i$ eigenvalue either. 

Now in dimension $\geq6$ it is known that $R$ has three irreducible parts: the scalar curvature, the traceless Ricci tensor and the Weyl tensor. We conclude the latter is 0, because the former are symmetric and hence cannot give a $4i$-eigenvalue. Finally, we recall the equivalence between Weyl and conformal flatness.

\subsubsection{The semi-Riemannian case}

Now suppose $(M,g)$ is an oriented $2n$-manifold and $g$ is an indefinite metric of signature $(2p,2q),\ \, p+q=n$. Let us denote 
\[  I_{p,q}=\left[ \begin{array}{cc} I_{2p} & 0 \\  0 &-I_{2q}   \end{array}\right]
\quad\mbox{and}\quad J_{p,q}=\left[ \begin{array}{cc} J_p & 0 \\  0 & -J_q \end{array}\right],
\quad\mbox{where}\quad J_p=\left[ \begin{array}{cc} 0& -I_p \\  I_p & 0 \end{array}\right]  .  \]
Thus each tangent space of $M$ admits an oriented orthonormal basis in which the metric is given by $I_{p,q}$. Next we consider the space  $F_{p,q}=SO_{2p,2q}/U_{p,q}$ whose elements are the linear complex structures compatible with the orientation and metric of semi-Euclidian space, or \textit{orthogonal} linear complex structures.
\begin{prop}\label{fpq_is_pK}
$F_{p,q}$ is a pseudo-K\"ahler symmetric space.
\end{prop}
\begin{proof}
Notice $\inv{J_{p,q}}={J_{p,q}}^T=-J_{p,q}$. $F_{p,q}$ is again a complex symmetric subspace of $GL_{2n}(\R)/GL_n(\C)$, because it is induced by the involutive automorphism $k\mapsto J_{p,q}k{J_{p,q}}^T$ of $SO_{2p,2q}$ with $U_{p,q}$ as the subgroup of fixed elements (we refer to the theory in \cite{KobNomi}). Since we have $T_JF_{p,q}$ identified with
\[  \m_J=\{ A\in \sol_{2p,2q}:\ AJ=-JA\}  \]
and the invariant complex structure is left multiplication by $J$, we have to check $JA\in\m_J$.
We know $AI_{p,q}=-I_{p,q}A^T$ and $JI_{p,q}J^T=I_{p,q}$. Hence
\[  JAI_{p,q}=-JI_{p,q}A^T=I_{p,q}J^TA^T=I_{p,q}(AJ)^T=-I_{p,q}(JA)^T   \]
as we wished. Since $kJA\inv{k}=kJ\inv{k}kA\inv{k}$, we have indeed an invariant complex structure.

Clearly $[[\m_J,\m_J],\m_J]\subset\m_J$, which is the condition for $\m_J$ to correspond to the canonical connection: a torsion free connection with parallel curvature. This is, moreover, the connection of the $SO_{2p,2q}$-invariant metric induced by the Killing form of $\sol_{2p,2q}$. Finally, if $\omega$ is the non-degenerate invariant pseudo-K\"ahler form, then
\[ \dx\omega(X,Y,Z)=\dx\omega(JX,JY,JZ)=\dx\omega(J^2X,J^2Y,J^2Z)=0  \]
and we are finished with the proof. Notice the curvature $R(A,B)C=-[[A,B],C]$. 
\end{proof}
Now we can talk about a new twistor space of $M$, also denoted ${\cal J}_+(M,g)=\{j\in{\cal J}(M)|\ j^*g=g\ \mbox{and $j$ induces the same orientation}\}$, with fibre $F_{p,q}$. We can also say it is the space of linear complex structures for which $g$ becomes type (1,1), or equivalently $g(j^+X,Y)=g(X,j^-Y)$.

By the remarks in the previous section, the equations of integrability of the almost complex structure $\jnab$ are the ones from theorem \ref{teo2.1} and precisely the same arguments from the definite case apply here.
\begin{teo}\label{teo2.2}
The twistor space ${\cal J}_+(M,g)$ is a complex manifold if, and only if, the metric is self-dual in case $2n=4$, or the metric is conformally flat in case $2n>4$.
\end{teo}
\begin{proof}
The semi-Riemannian decomposition of the curvature tensor is sustained in all signatures and, according to \cite{Besse}, theorem 1.165, the vanishing of the semi-Riemannian Weyl tensor corresponds to conformal flatness. In dimension $4$, the case for $SO_{2,2}$ also resumes to self-duality ($W^-=0$) because the Hodge operator still verifies $*^2=1$ and this group is not simple.
\end{proof}

\subsection{Holomorphic maps into twistor space}

Let $\cal Z$ be any of the previously described twistor spaces over a manifold $(M,\nabla)$. Suppose $(N,J^N)$ is a given almost complex manifold and $\psi:N\rr\cal Z$ a given map. Let $f=\pi\circ\psi$ and let $\psi^*\Phi$ be the pullback of the tautological almost complex structure of the bundle $E$ described in (\ref{seqexacta}): $\psi^*\Phi_x$ agrees with $\psi(x)$ for all $x\in N$. This induces a decomposition $\psi^*E\otimes\C=\psi^+\oplus\psi^-$ into $\pm i$-eigenbundles. Now we need a lemma whose proof was already given in two particular situations: in \cite{Rawnsley} for the Riemannian case and in \cite{Albu} for the symplectic case. It is a result of a technical sort, which carries straightforwardly to the present setting.
\begin{lema}[\cite{Rawnsley}]\label{lematech}
On any twistor space the following conditions are equivalent:\\
(i) \,$\psi$ is $(J^N,\jnab)$ pseudo-holomorphic.\\
(ii) \,$\dx f\circ J^N=\psi^*\Phi\circ\dx f$\, and \,$(f^*\nabla\,_u\psi^*\Phi)(\psi^+)=0,\ \,\forall u\in T^+N$.\\
(iii) \,$\dx f(T^+N)\subset\psi^+$\, and \,$f^*\nabla\,_u(\Gamma\psi^+)\subset(\Gamma\psi^+),\ \,\forall u\in T^+N$.
\end{lema}
Now suppose $N=M$ and $\psi=J:M\rr \cal Z$ is a smooth section. Let $J$ itself play the role of $J^N$ above, as it is an almost complex structure on $M$. Then $f=\Id$ and $J^*\Phi=J$. Moreover, the space of sections $\Gamma J^+=\Gamma T^+M=\XIS^+$.

The following result generalizes one from \cite{Sal} in two ways.
\begin{prop}    \label{prop_com_sal}
For $\nabla$ torsion free, the almost complex structure $J$ is integrable if the map $J$ is $(J,\jnab)$ pseudo-holomorphic.

For the semi-Riemannian twistor space with the Levi-Civita connection, the condition is also sufficient.
\end{prop}
\begin{proof}
Let us analyse (iii) in the lemma. The first part holds trivially and the second resumes to
\begin{equation}\label{condition_B}
\nab{u}{v}\in\XIS^+,\ \ \forall u,v\in\XIS^+.
\end{equation}
But then the integrability follows by the vanishing of the Nijenhuis tensor, which is well known to be equivalent to
\[  [u,v]=\nab{u}{v}-\nab{v}{u}\in \XIS^+.  \]

Now suppose we are in the semi-Riemannian setting and the last equation is fullfield, ie. $[\XIS^+,\XIS^+]\subset \XIS^+$, which is the same as $J$ being integrable. By hypothesis the metric $g$ is type (1,1) relatively to $J$. Let us define a 3-tensor $\Theta(u,v,w)=g(\nab{u}{v},w)$ in $\XIS^+$. It is indeed $\cinf{M}(\C)$-linear in $v$ because $g(v,w)=0$. By the same reason and the fact that $\nab{}{g}=0$, $\Theta$ is skew-symmetric in $v,w$:
\[   g(\nab{u}{v},w)=u\cdot g(v,w)-g(v,\nab{u}{w})=-g(v,\nab{u}{w}).  \]
But the integrability of $J$ implies $\Theta$ is symmetric in $u,v$. These two conclusions lead to $\Theta=0$ and therefore (\ref{condition_B}) is valid again. Applying the lemma, we see $J$ is pseudo-holomorphic.
\end{proof}

\subsection{Affine transformations of twistor space}

Let $M,M_1$ be two manifolds and $\sigma:M\rightarrow M_1$ a diffeomorphism. Then $\sigma$ induces an invertible transformation from ${\cal J}(M)$ onto ${\cal J}(M_1)$ preserving the fibres, ie. a map $\Sigma$ such that the diagram
\begin{eqnarray*}
\begin{array}{ccc}
{\cal J}(M) & \stackrel{\Sigma}\longrightarrow & {\cal J}(M_1) \\
\pi\downarrow \hspace{3mm}  & & \hspace{3mm} \downarrow\pi_1 \\
M & \stackrel{\sigma}\longrightarrow & M_1
\end{array}
\end{eqnarray*}
commutes. Indeed, for any $y\in M_1,\ j\in\pi^{-1}(\sigma^{-1}(y))$ we define
\begin{equation}\label{definideSigma}
 \Sigma(j)=\dx\sigma\circ j\circ\dx\sigma^{-1}
\end{equation}
which is an element in $\pi_1^{-1}(y)$. It is trivial to check $\Sigma$ is well defined. 

We may suppose furthermore that $\sigma$ preserves some extra $G$-structure, in the sense that it interchanges the principal $G$-bundle of frames of $M$ and $M_1$. Then it induces a map $\Sigma:{\cal Z}\rr{\cal Z}_1$ between the twistor subspaces whose fibres are $G/G\cap GL_n(\C)$.

\newcommand{\jnabum}{{\cal J}^{\nabla^1}}       
\newcommand{\jnabdois}{{\cal J}^{\nabla^2}}       

Assume we have twistor almost complex structures $\jnab$ and $\jnabum$, on the respective twistor spaces, where $\nabla^1=\sigma\cdot\nabla$ and $\nabla$ is any given linear $G$-connection on $M$. Recall that for any $Z,W$ vector fields on $M_1$,
\[ (\sigma\cdot\nabla)_ZW=\sigma\cdot\bigl(\nab{\inv{\sigma}\cdot Z}{\inv{\sigma}\cdot W}\bigr) \]
where $\sigma\cdot X\,_y=\dx\sigma(X_{\inv{\sigma}(y)}),\ \,\forall y\in M_1$. The new connection is again a linear $G$-connection, and $\sigma$ becomes an affine transformation. Since one can also see $\Sigma$ as the map $\sigma\cdot\,$ acting on twistors, the following must be true.
\begin{teo}[\cite{AR}]  \label{teo4.1}
$\Sigma:{\cal Z}\rr{\cal Z}_1$ is pseudo-holomorphic.
\end{teo}
\begin{proof}
This proof is considerably shorter than the one in the reference. Notice that $\Sigma$, when restricted to each fibre, extends to a linear map between
$\End{T_{\sigma^{-1}(y)}M}$ and $\End{T_yM_1}$. Hence applying (\ref{estrcomplexatwistor})
\[ \dx\Sigma(jA)=\Sigma(jA)=\Sigma(j)\Sigma(A)=\Sigma(j)\,\dx\Sigma(A)   \]
and we may conclude the map is {\it vertically} pseudo-holomorphic.

Now we shall check part (ii) of lemma (\ref{lematech}) considering $\Sigma$ as a map into the second twistor space ${\cal Z}_1$. Let $f=\sigma\circ\pi=\pi_1\circ\Sigma$. By definition, for any $X\in T_j{\cal J}(M)$ we have
\[  \dx f\circ\jnab (X)=\dx\sigma\circ\dx\pi(\jnab X)=\dx\sigma\circ j\circ(\dx\sigma^{-1}\dx\sigma)\circ\dx\pi X=\Sigma(j)\dx f(X)   \]
which is the first part of the condition. For the second we take $u\in{\hnab}^+,\ \Phi,\Phi^1$ the canonical sections (cf. proposition \ref{horizontais}) and notice
\begin{equation}\label{eeeee}
(f^*\nabla^1_u\Sigma^*\Phi^1)\Sigma^+= ((\Sigma^*\pi_1^*\nabla^1)_u\Sigma^*\Phi^1)\Sigma^+ = ((\pi_1^*\nabla^1)_{\Sigma_*u}\Phi^1)\Sigma^+  
\end{equation}
so the theorem follows after the proof that $\Sigma_*\hnab={\cal H}^{\nabla^1}$.
This turns out to be exactly the case when we consider the particular connection $\nabla^1$.

Notice that $\Sigma^*\Phi^1_j=\Phi^1_{\Sigma(j)}=\dx\sigma j\inv{\dx\sigma}=\sigma\cdot\Phi_j$. Also it is not difficult to compute the formula, for any section $\xi$ of $\sigma^*TM_1$,
\[   \sigma^*\nabla^1_Z\xi=\sigma^*\bigl(\sigma\cdot(\nab{Z}{\inv{\sigma}\cdot \xi)}\bigr)  \]
for any $Z\in TM$. Finally suppose $X\in\hnab$. According to proposition (\ref{horizontais}) we have $\pi^*\nab{X}{\Phi}=0$ and want to prove a similar equality for $\Sigma_*X$. Now
\begin{eqnarray*}
\pi_1^*\nabla^1_{\Sigma_*X}{\Phi^1} & = & \bigl((\pi_1\circ\Sigma)^*\nabla^1\bigr)_X\Sigma^*\Phi^1\ =\ \bigl((\sigma\circ\pi)^*\nabla^1\bigr)_X\Sigma^*\Phi^1   \\ 
      & =  &  (\pi^*\sigma^*\nabla^1)_X\sigma\cdot\Phi  \ = \
   \pi^*\sigma^*\bigl(\sigma\cdot(\pi^*\nab{X}{(\inv{\sigma}\cdot\sigma\cdot\Phi)})\bigr) \ =\ 0
\end{eqnarray*}
as we wished.
\end{proof}
The principle behind the last computation is the fact that an affine transformation sends $\nabla$-horizontal frames into $\nabla^1$-horizontal frames. Now suppose we have on $M_1$ a second linear connection $\nabla^2=\nabla^1+{\cal A}$.
\begin{coro}
The map $\Sigma:({\cal Z},\jnab)\rr({\cal Z}_1,\jnabdois)$ is pseudo-holomorphic if, and only if, $j_1^+{\cal A}_{j_1^-Y}j_1^-=0$, \,$\forall Y\in TM_1,\ \forall j_1\in {\cal Z}_1$.
\end{coro}
\begin{proof}
We know that for any $u\in{\hnab_j}^+$, such that $j\in{\cal Z}$, we have $\Sigma_*u=v\in{{\cal H}^{\nabla^1}_{\Sigma(j)}}^+$. So we just have to follow the last proof from that point of formula (\ref{eeeee}), which must vanish:
\begin{eqnarray*}
((\pi_1^*\nabla^2)_{\Sigma_*u}\Phi^1)\Sigma^+=0 \qquad\qquad\qquad\qquad \\
\Longleftrightarrow\ \ [\pi_1^*{\cal A}_{v},\Phi^1]\Sigma^+=0 \ \ 
\Longleftrightarrow\ \ [{\cal A}_{\dx{\pi_1}_{j_1}(v)},j_1]j_1^+=0, \ \,\forall j_1\in{\cal Z}_1.
\end{eqnarray*}
By definition $\dx{\pi_1}_{j_1}(v)=Y-ij_1Y\in T_{\pi_1(j_1)}^+M_1$ for some $Y\in TM_1$. Since
\begin{eqnarray*}
[{\cal A},j_1]j_1^+ &  = & \bigl((j_1^++j_1^-){\cal A}j_1-j_1(j_1^++j_1^-){\cal A}\bigr)j_1^+  \\ & = & i\bigl(j_1^+{\cal A}j_1^++j_1^-{\cal A}j_1^+-j_1^+{\cal A}j_1^++j_1^-{\cal A}j_1^+\bigr) \ =\ 2i\,j_1^-{\cal A}j_1^+
\end{eqnarray*}
the condition on ${\cal A}$ is equivalent to $j_1^-{\cal A}_{j_1^+Y}j_1^+=j_1^+{\cal A}_{j_1^-Y}j_1^-=0$.
\end{proof}
Notice that if $\sigma=\Id$, then $\Sigma=\Id$; hence the corollary gives the necessary and sufficient condition on ${\cal A}$ in order to have $\jnab=\jnabdois$. From this remark one proves easily that the twistor almost complex structure on the semi-Riemannian twistor space is independent of a conformal change of the metric, a well known result in the definite case (\cite{Obri}). Just recall the difference tensor ${\cal A}=\nabla^2-\nabla$ induced by the metrics $g$ and $e^{2f}g$ is given by ${\cal A}_XY=X(f)Y+Y(f)X-g(X,Y)\mathrm{grad}\,f$.

Also we remark that theorem (\ref{teo4.1}) is coherent with the integrability equations of (\ref{teo2.1}) because $\Sigma(j)^\pm=\Sigma(j^\pm),\ \,\forall j$, and the torsion and curvature tensors satisfy $T^{\sigma\cdot\nabla}=\sigma\cdot T$ and $R^{\sigma\cdot\nabla}=\sigma\cdot R$.

\begin{coro}
Suppose $\sigma$ is an isometry of a semi-Riemannian manifold $(M,g)$. Then the map $\Sigma:{\cal J}_+(M,g)\rr{\cal J}_+(M,g)$ is pseudo-holomorphic.
\end{coro}
\begin{proof}
The affinely transformed connection $\sigma\cdot\nabla$ of the Levi-Civita connection $\nabla$ is also a metric and torsion free connection. By uniqueness, the two connections coincide.
\end{proof}

\section{The case for the pseudo-spheres}

\subsection{Useful results}

Now we consider the $2n$-dimensional pseudo-sphere $S^{2n}_{2q}=SO_{2p+1,2q}/SO_{2p,2q}$ with its usual $SO_{2p+1,2q}$-invariant metric $\left\langle\ ,\ \right\rangle$, where $n=p+q,\ p,q\geq 0$. We concede to the usual prefix `pseudo', remarking it is not refering to the complex manifold terminology. Notice the invariant metric induced by the Killing form is the same as the metric of the flat semi-Euclidian space $\R^{2p+1,2q}$ restricted to the tangent bundle of the homogeneous space of norm 1 vectors. Also recall that this even dimensional pseudo-sphere is diffeomorphic to $S^{2p}\times\R^{2q}$. We let $\zo^{p,q}$ denote the twistor space ${\cal J}_+(S^{2n}_{2q}, \left\langle\ ,\ \right\rangle)$.

Recall $S^{2n}_{2q}$ is a connected, simply-connected complete semi-Riemannian manifold of constant sectional curvature 1. Hence all twistor spaces $\zo^{p,q}$ are complex manifolds. 
\begin{prop}\label{naoe'psdKahler}
$S^{2n}_{2q}$ cannot be a pseudo-K\"ahler manifold for any complex structure compatible with the metric, except if $p+q=1$.
\end{prop}
\begin{proof}
Let $q=0$ and $p>1$. Then the Riemannian spheres are not K\"ahler by topological reasons (a closed K\"ahler form yields a manifold with no volume). 

Now suppose both $p,q>0$. Then $S^{2n}_{2q}$ cannot be pseudo-K\"ahler because of the classification of space-forms of this kind. Consider the open subset $\C\Pro^n_q$ of complex projective space consisting of lines generated by $z\in\C^{n+1}$ such that 
\[  \sum_{i=0}^pz_i\overline{z}_i-\sum_{i=p+1}^nz_i\overline{z}_i  \]
is greater than 0. Then, for any $c>0$, this space inherits an indefinite K\"ahler metric of constant holomorphic sectional curvature $c$. Now a result of \cite{BR} says that a connected, simply-connected, complete pseudo-K\"ahler manifold of signature $(2p,2q)$ and constant holomorphic sectional curvature $c$ must be isometric and biholomorphic to $\C\Pro^n_q$. So the pseudo-sphere should be isometric to this projective subspace, with $c=1$, because its sectional, and hence holomorphic sectional, curvature is constant 1. However, this is in contradiction with the fact that not all the sectional curvatures of $\C\Pro^n_q$ are 1. Indeed for any $X,Y$ tangent to this manifold, with $\left\langle X,X\right\rangle=1,\ \left\langle Y,Y\right\rangle =-1$ and $\left\langle X,Y\right\rangle =0$, then 
\[  R(X,JX,X,JX)=1\qquad\mbox{and}\qquad R(X,Y,X,Y)=-\dfrac{1}{4} , \]
as we can see by a formula of \cite{BR}. One may also argue that the two spaces are in fact not homotopically equivalent if $p>1$.
\end{proof}
The twistor spaces of pseudo-spheres are not very difficult to describe. 
\begin{teo}
The following are biholomorphic identities:
\[   \zo^{p,q}=\frac{SO_{2p+1,2q}}{U_{p,q}}=\frac{SO_{2p+2,2q}}{U_{p+1,q}} .  \]
\end{teo}
\begin{proof}
By theorem \ref{teo4.1} the Lie group $SO_{2p+1,2q}$ acts by biholomorphisms on $\zo^{p,q}$. The isotropy subgroup is evidently $U_{p,q}$ as we deduce from the definition (\ref{definideSigma}). By counting dimensions, the first identity follows. We note that this action can be seen, locally, as $b\cdot(x,j)=(bx,bj\inv{b})\in SO_{2p+1,2q}/SO_{2p,2q}\times SO_{2p,2q}/U_{p,q}$. 

For the second identity, we note that every $j\in \inv{\pi}(x)\subset\zo^{p,q}$ extends to a linear complex structure in $\R^{2p+2,2q}=\R1+\R^{2p+1,2q}$, writing $\overline{j}(x)=-1,\ \overline{j}(1)=x$. This extension is in fact the \textit{identity} map, since for any linear orthogonal complex structure $J$ in $\R^{2p+2,2q}$ we get 
\[ \left\langle 1,J(1)\right\rangle=-\left\langle J(1),1\right\rangle =0   \] 
and due to the conjugation of $J$ by a $b\in SO_{2p+1,2q}$ agreeing with the action above. Notice the bundle projection to the pseudo-sphere is $J\mapsto J(1)$.
\end{proof}
Here is a well known result whose proof, at the light of the theorem, might be interesting to notice.
\begin{coro}
$\C\Pro^3$ is the twistor space of the 4-sphere.
\end{coro}
\begin{proof}
We recall the Riemannian twistor bundle is usually seen as $\Hamil^2/\C^*\rr\Hamil\Pro^1=S^4$ so the whole space is $\C\Pro^3$ and the fibre is $\C\Pro^1$. The latter agrees with the 2-sphere of normed 1, self dual 2-forms. Now the \textit{holomorphic} identification of 3-projective space with $SO_6/U_3$ comes from a special isomorphism $\s\uni(4)\simeq\sol(6)$ \,(cf. \cite{Hel}, pp. 518-519, the coincidence AIII(p=3,q=1)=DIII(n=3)).
\end{proof}

It is known by a result of A.~Borel and J. P.~Serre that the only spheres which admit almost complex structures are $S^2$ and $S^6$. The results presented above lead to a new proof of the following interesting result of C. Lebrun.
\begin{teo}[\cite{Lebr}]
There is no integrable orthogonal complex structure on $S^6$.
\end{teo}
\begin{proof}
Suppose there exists a section $J:S^6\rr\zo^{3,0}$ representing such an integrable complex structure. By the existence of local complex charts, $J$ must me a smooth section. It is also holomorphic by proposition \ref{prop_com_sal}. Thus $S^6$ embbeds as a complex submanifold of the K\"ahler manifold $SO_{8}/U_{4}$, and hence it is itself a K\"ahler manifold --- a contradiction.
\end{proof}

\subsection{The metric on $\zo^{p,q}$}

The spaces ${\cal J}_+(M,g)$ inherit a metric $a\pi^*g+bg_f$, where $g_f$ is the invariant metric defined on the fibres via the connection and $a,b$ are any two non-vanishing functions. This works for any manifold and yields a metric compatible with $\jnab$, as it is simple to check.

In the present application to pseudo-spheres we shall find $a,b$ such that the metric on $\zo^{p,q}=F_{p+1,q}$ agrees with the $SO_{2p+2,2q}$-invariant one of proposition \ref{fpq_is_pK}. Let $1$ represent a norm 1 direction in semi-Euclidian space and let
$\m^p_J=\{ A\in \sol_{2p,2q}:\ AJ=-JA\}$. Since the bundle projection is given by the \textit{linear} map $\pi(J)=J(1)$, it is easy to see that the vectors tangent to the fibres, ie. those in ${\cal V}=\ker\dx\pi$, correspond to
\[   A\in \m^{p+1}_J\qquad\mbox{such that}\qquad A1=AJ1=0.   \]
It follows that, for any $X\in T_{J(1)}S^{2n}_{2q}$, we get $\left\langle AX,1\right\rangle =\left\langle AX,J1\right\rangle =0$. Hence, a tangent vector $A\in\m^{p+1}_J$ is tangent to the fibres of the twistor bundle if $A$ coincides with an endomorphism of $\{1,J1\}^\perp$. We shall denote the vertical part of any tangent $A$ by $A'$.
\begin{lema}
The Killing form of $\sol_{k,l}$ is given by $B_{k,l}(A_1,A_2)=(k+l-2)\Tr{A_1A_2}$.
\end{lema}
\begin{proof}
It is well known the Killing form of $\sol(k+l,\C)=\g$ is given by the formula above. On the other hand, for any real form $\g_0$ of a complex Lie algebra, ie. any real Lie algebra such that $\g_0\otimes\C=\g$, its Killing form is clearly the restriction to real vectors of the Killing form of $\g$. So we just have to prove $\sol_{k,l}$ is a real form of $\g$. Given $X_1\in\sol_k,\ X_2\ \mbox{any}\ k\times l\ \mbox{matrix, and}\ X_3\in\sol_l$, the map
\[\left[ \begin{array}{cc} X_1 & X_2 \\ X_2^T & X_3 \end{array}\right] 
\mapsto\left[ \begin{array}{cc} X_1 & iX_2 \\ -iX_2^T & X_3 \end{array}\right]   \] 
can easily be seen to be an isomorphism of Lie algebras. Of course its image is a real form of $\sol(k+l,\C)$, and since isomorphisms induce isometries for the Killing metric, we are finished (cf. \cite{Hel}, pp 189, 239 for details).
\end{proof}
Returning to the above, we write $\left\langle A_1,A_2\right\rangle_k=-B_{2k,2l}(A_1,A_2)$ (recall the Killing form is negative definite on the compact orthogonal Lie algebra). Now computing the trace using a basis containing $1$ and $J1$, we find
\begin{equation}
(2p+2q)\left\langle A_1',A_2'\right\rangle_p =(2p+2q-2)\left\langle A_1',A_2'\right\rangle_{p+1} 
\end{equation}
for any vertical vectors $A_1',A_2'$. We have proved part of the following result.
\begin{prop}
For any vectors $A,B\in T_JF_{p+1,q}=\m_J^{p+1}$, we have
\begin{equation}\label{metricahv}
\left\langle A,B\right\rangle_{p+1}=8n \left\langle A1,B1 \right\rangle +\frac{n}{n-1}\left\langle A',B'\right\rangle_p .
\end{equation}
In particular, the index $i_{p,q}$ of the metric on $F_{p,q}$ (the number of time-like vectors in an orthonormal basis) is $q^2-q+2pq$.
\end{prop}
\begin{proof}
Let $\{X_1,\ldots,X_n,JX_1,\ldots,JX_n\}$ be a (direct) orthonormal basis of $\{1,J1\}^\perp$ in $\R^{2p+2,2q}=\R1+\R^{2p+1,2q}$, let $\epsilon_i=\left\langle X_i,X_i\right\rangle =\left\langle JX_i,JX_i\right\rangle$ and set, for $1\leq i\leq n$,
\[  A_i1=\epsilon_i X_i,\qquad A_iJ1=-JA_i1=-\epsilon_i JX_i,\qquad A_iX_j=-\delta_{ij}1,\qquad A_iJX_j=\delta_{ij}J1  \]
for any $j\leq n$. Then clearly $A_iJ=-JA_i$, and $A_i\in\sol_{2p+2,2q}$ because
\[  \left\langle A_i1,X_j\right\rangle =\epsilon_i\left\langle X_i,X_j\right\rangle=\delta_{ij}=-\left\langle 1,A_iX_j\right\rangle , \]
\[  \left\langle A_i1,JX_j\right\rangle =\epsilon_i\left\langle X_i,JX_j\right\rangle=0=-\left\langle 1,A_iJX_j\right\rangle. \]
Also $\left\langle A_i1,1\right\rangle =\epsilon_i\left\langle X_i,1\right\rangle=0=-\left\langle 1,A_i1\right\rangle $ with equal conclusion for $J1$. Finally $\left\langle A_iX_k,X_j\right\rangle =0=-\left\langle X_k,A_iX_j\right\rangle$ as we wished.

It is clear enough that $A_i'=0$. Now we extend the set of endomorphisms $A_1,\ldots,A_n$ to a basis of the horizontal tangent bundle $\hnab$ putting $A_{i+n}=JA_i$.

If we compute the horizontal part $\left\langle \pi_*A_i,\pi_*A_j\right\rangle$ of the metric, we get
\[   \left\langle A_i1,A_j1 \right\rangle=\epsilon_i\delta_{ij} .  \]
On the other hand, computing directly $\left\langle A_i,A_j\right\rangle_{p+1}$ we get, for $i,j\leq n$,
\begin{eqnarray*}
-(2n+2-2)\Tr{A_iA_j} & = & -2n\bigl(\left\langle A_iA_j1,1\right\rangle+\left\langle A_iA_jJ1,J1\right\rangle + \\  & & \ \ \ +\sum_{k=1}^n(\epsilon_k\left\langle A_iA_jX_k,X_k\right\rangle+\epsilon_k\left\langle A_iA_jJX_k,JX_k\right\rangle)\bigr) \\
& = & +2n\bigl(2\left\langle X_j,X_i\right\rangle\epsilon_i\epsilon_j+\sum_k2\epsilon_k\left\langle A_jX_k,A_iX_k\right\rangle\bigr)  \\
& = & 4n (\epsilon_i\delta_{ij}+\sum_k\epsilon_k\delta_{kj}\delta_{ki})=8n\,\epsilon_i\delta_{ij} 
\end{eqnarray*}
which leads to formula (\ref{metricahv}). It is easy to prove $\Tr{JA_iA_j}=0$ using the same basis, and clearly $\Tr{JA_iJA_j}=\Tr{A_iA_j}$. Also worth noticing is that $\Tr{A_iA'}=0$ for any vertical vector $A'$. The formula for the index follows by induction; we have $i_{0,q}=q(2q-1)-q^2=q^2-q$ and $i_{p+1,q}=i_{p,q}+2q$, therefore $i_{p,q}=q^2-q+2pq$.
\end{proof}

\subsection{Old and new formulas for $\dx\omega$}

Suppose $(M,g)$ is a semi-Riemannian manifold, $\nabla$ is the Levi-Civita connection and $\cal Z$ is its twistor space. Let $\ ':T{\cal Z}\rr{\cal V}\subset\End{E}$ be the projection with kernel $\hnab\simeq E$. Then this projection can be seen as a 1-form on $\cal Z$ and thus capable of inducing a translation of the usual connection in $E$ to a pseudo-unitary connection:
\begin{equation}
D_A=\pi^*\nab{A}{\,}-A'.
\end{equation}
Indeed, since $g(A'X,Y)=-g(X,A'Y)$, for all $X,Y\in TM$, $D$ is a metric connection for the natural metric $\pi^*g$ in $E$, and from proposition \ref{horizontais} it follows that $D\Phi=0$. Moreover, $D$ preserves $\cal V$ and therefore we find, as in \cite{Obri}, a new linear connection, also denoted by $D$, on the tangent bundle of $\cal Z$ preserving the decomposition $\hnab\oplus\cal V$. Still, $D\jnab=0$.

It is known that the torsion
\[  T^D(A,B)=\pi^*T^\nabla_{\,A,B}-A'\pi_*B+B'\pi_*A+(\pi^*R^\nabla_{\,A,B})' \]
--- this was computed in the general setting in \cite{Obri} and of course holds in the present case (for which $\pi^* T^\nabla=0$). Notice also the horizontal and vertical parts decomposition. Furthermore, the formula leads to a proof of theorem \ref{teo2.1} which we succintely recall: using a well known identity for the Nijenhuis tensor, $N(A,B)=8\,\prl{\jnab}^+T^D({\jnab}^-A,{\jnab}^-B)$ for a complex connection, equations (\ref{int}) follow with little extra work.

Now let
\begin{equation}\label{gt}
G^t(A,B)=8n\,\pi^*g(A,B)+t\,g_f(A',B')
\end{equation}
be the metric on the twistor bundle defined via the connection ($t\in\R\backslash\{0\}$). As we have seen, $g_f(A',B')$ essentially agrees with the trace ($(2n-2)$ times), so it is simple to verify $Dg_f$, and hence $DG^t$, is zero. We may also define a non-degenerate parallel 2-form $\Omega=G^t(\jnab\ ,\ )$.
\begin{prop}\label{dOmega1}
Let $A,B,C\in\hnab\cup\cal V$. Then $\dx\Omega(A,B,C)\neq0$ if, and only if, two of the vectors are horizontal and the other is vertical. If $X,Y\in\hnab_j,\ A\in{\cal V}_j$, then 
\begin{equation}\label{dOmega}
  \dx\Omega_j(X,Y,A)=-16n\, g(jAX,Y)+t\, g_f(jR^\nabla_{\,X,Y},A)
\end{equation}
where we identify $X$ with $\pi_*X\in T_{\pi(j)}M$.
\end{prop}
\begin{proof}
It is known that, for connections such that $D\Omega=0$, we have
\[     \dx\Omega(A,B,C)=\cyclic_{A,B,C}\Omega(T^D(A,B),C) . \]
Hence the result follows by carefull thinking of all four cases of horizontal and vertical choices.
Therefore (\ref{dOmega}) is deduced from
\begin{eqnarray*}
\Omega_j(T^D(X,Y),A)+\Omega_j(T^D(A,X),Y)+\Omega_j(T^D(Y,A),X) \qquad\\
 = t\,g_f(jR^\nabla_{\,X,Y},A)-8n\,g(jAX,Y)+8n\,g(jAY,X) 
\end{eqnarray*}
which is the same as above. It is important to notice we are only using the vertical part of $R^\nabla_{\,X,Y}$, ie. the one which anti-commutes with $J$, by the reason that it is perpendicular to $\uni_{p,q}$ with respect to the trace.
\end{proof}
Let $J:M\rr\cal Z$ be a smooth section and let $\omega$ denote the associated 2-form $g(J\ ,\ )$. Then $J^*\Omega=8n\,\omega+J^*\tau$ where $\tau$ denotes the vertical part.
\begin{prop}\label{domega1}
Suppose $\dx\Omega=0$. Then
\begin{equation}\label{domega}
\dx\omega(X,Y,Z)=\cyclic_{X,Y,Z}\,-\frac{t}{16n}\,g_f(R^\nabla_{\,X,Y},\nab{Z}{J}) .
\end{equation}
\end{prop}
\begin{proof}
As we have seen earlier, in section 2.1, the vertical part of $\dx J(X)$ is $\frac{1}{2}J\nab{X}{J}$. The computations above also show what the result of $\dx \omega=-\frac{1}{8n}J^*\dx\tau$ must be: for all $X,Y,Z\in TM$, we have $\dx\omega(X,Y,Z)$ equal to
\begin{eqnarray*}
 -\frac{1}{8n}\,\dx\tau_J(J_*X,J_*Y,J_*Z) 
& = & -\frac{1}{8n}\cyclic_{X,Y,Z}\,\tau_J(T^D(J_*X,J_*Y),J_*Z)  \\ 
& = & -\frac{t}{16n}\cyclic_{X,Y,Z}\,g_f(\jnab_J\pi^*R^\nabla_{J_*X,J_*Y},J\nab{Z}{J})
\end{eqnarray*}
and the result follows.
\end{proof}
Notice we can consider a 2-form on the twistor space $\varpi=\pi^*g(\jnab\ ,\ )$ and the pull-back of this by $J$ agrees with $\omega$. Then it is not hard to see, as in proposition \ref{dOmega1}, that $J^*\dx\varpi$ leads to the \textit{old} formula
\begin{equation}\label{formula simples e estranha}
   \dx\omega(X,Y,Z)=\cyclic_{X,Y,Z}g((\nab{Z}{J})X,Y)
\end{equation}
which is not so easy to deduce if we apply directly the Levi-Civita connection.

We easily discover that $\dx\varpi$ depends on one vertical and two horizontal vector fields (cf. proposition \ref{dOmega1}). For instance,
\begin{eqnarray*}
\dx\varpi(B1,C1,A') & = & \varpi(T^D(B1,C1),A')+\varpi(T^D(C1,A'),B1)+\varpi(T^D(A',B1),C1)  \\
& = & \varpi(A'C1,B1)-\varpi(A'B1,C1)\ =\ -2\varpi(A'B1,C1).
\end{eqnarray*}
We show the following proposition in order to understand better this 3-form.
\begin{prop}
$\dx\varpi$ is a form of type (1,2)+(2,1).
\end{prop}
\begin{proof}
Suppose $X-ijX\in{\hnab}_j^+,\ A'-ij A'\in{\cal V}_j^+$. Then, in computing $\dx\varpi^{(3,0)}$ by the formula above, we would cross with the computation
\[  (A'-ijA')(X-ijX)=A'X-jA'jX-i(jA'X+A'jX)=0  \]
which yields the conclusion that part must vanish. If $\dx\varpi^{1,2}=0$, then we would have $\dx\varpi=0$ in contradiction with the above.
\end{proof}

\subsection{Application to the pseudo-spheres}

We return to the study of the bundle $\zo^{p,q}\rr S^{2n}_{2q}$. By the result of (\ref{metricahv}) in section 3.2 we have an $SO_{2p+2,2q}$-invariant metric compatible with the complex structure $\jnab$, which yields an identification $\zo^{p,q}=F_{p+1,q}$. We recall the decomposition of $A\in T\zo^{p,q}$ as
\[  A=A1+A'   \]
into horizontal and vertical directions. However, if we take coordinates $(x^1,\ldots,x^{2n})$ on $S^{2n}_{2q}$, then we still denote the horizontal vector field $\inv{(\dx\pi)}(\partial/\partial x^i)$ by $\partial_i$.

As explained in section 3.3 we may define a new linear connection $D$ on $\zo^{p,q}$, preserving the splitting $\hnab\oplus{\cal V}$. We start by checking the expression for the torsion in general terms, since the result in \cite{Obri} is capable of further improvement. The vertical part is\footnote{The reader must distinguish between Lie and commutator brackets.}
\begin{eqnarray*}
T^D(A,B)'& = & D_AB'-D_BA'-[A,B]'  \\
& = & \frac{1}{2}\biggl(D_A(\Phi\pi^*\nabla_B\Phi)-D_B(\Phi\pi^*\nabla_A\Phi)- \Phi\pi^*\nabla_{[A,B]}\Phi\biggr)  \\
& = & \frac{1}{2}\Phi\biggl( [R^{\pi^*\nabla}_{A,B},\Phi]-[A',\pi^*\nabla_B\Phi]+ [B',\pi^*\nabla_A\Phi]\biggr) \\
& = & \frac{1}{2}\Phi\bigl( -2\Phi(R^{\pi^*\nabla}_{A,B})'+2[A',\Phi B']-2[B',\Phi A']\bigr) \ =\ (R^{\pi^*\nabla}_{A,B})'.
\end{eqnarray*}
and the horizontal part is quickly checked for three cases: for two horizontal vectors, $T^D(\partial_i,\partial_j)=\pi^*T^\nabla(\partial_i,\partial_j)$, for two verticals we have
$T^D(A',B')1=0$ because the vertical tangent bundle $\cal V$ is integrable and $D$ preserves $\cal V$. Last, but not least,
\begin{eqnarray*}
T^D(A',\partial_i)1 & = & (\pi^*\nab{\partial_i}{A'}-\pi^*\nab{A'}{\partial_i}+A'-[A',\partial_i])1\\
& = & -A'\pi^*\nabla_{\partial_i}{1}-\pi^*\nab{A'}{\partial_i}-A'\partial_i\ =\ -A'\partial_i
\end{eqnarray*}
and thus, in sum, $T^D(A,B)1=-A'B1+B'A1$.

\subsubsection{Non-existence of orthogonal complex structures}

Now suppose $J:S^{2n}_{2q}\rr\zo^{p,q}$ is an integrable complex structure and let $\omega$ denote the associated 2-form. Then $\dx\omega$ is type (1,2)+(2,1) because $\omega$ is type $(1,1)$ and because $d=\partial+\db$. Also, recall $\dx J$ preserves types by proposition \ref{prop_com_sal}. We are going to use the formula (\ref{domega}) with $R_{X,Y}Z=\left\langle Y,Z\right\rangle X-\left\langle X,Z\right\rangle Y$. We therefore must check carefully the weights of the metric. We saw in (\ref{metricahv}) that the pseudo-K\"ahler metric of the twistor space is the metric $G^t$ from (\ref{gt}) with
\[  t=\frac{n}{n-1}.   \]
Since $g_f$ on the fibre is $-(2n-2)\Tr{}$, we find by proposition \ref{domega1}
\[  \dx\omega(X,Y,Z)=\cyclic_{X,Y,Z}\,\frac{1}{8}\,\Tr{(R^\nabla_{\,X,Y}\nab{Z}{J})} .  \]
\begin{prop}
$\dx\omega=0$.
\end{prop}
\begin{proof}
Let $\{X_1,\ldots,X_n,JX_1,\ldots,JX_n\}$ be a local orthonormal frame of the tangent bundle $\{1,J1\}^\perp$ in $\R^{2p+2,2q}=\R1+\R^{2p+1,2q}$ and let $\epsilon_k=\left\langle X_k,X_k\right\rangle =\left\langle JX_k,JX_k\right\rangle$. For any real endomorphism $C$ of $T_xS^{2n}_{2q}$ we have
\[ {\mathrm{Tr}}_{\R}\,C\ =\ \prl\epsilon_k\langle Ce_k,\overline{e_k}\rangle\ =\ \epsilon_k\langle CX_k,X_k\rangle+\epsilon_k\langle CJX_k,JX_k\rangle  \]
where $e_k=X_k-iJX_k$ (repeated indices represent a sum from 1 to $n$). Notice $\langle e_k,\overline{e_k}\rangle=2\epsilon_k$. Hence
\begin{eqnarray*}
\dx\omega(X,Y,Z) & = & \cyclic_{X,Y,Z}\,\frac{1}{8}\,\prl\epsilon_k\langle R^\nabla_{X,Y}(\nab{Z}{J})e_k,\overline{e_k}\rangle\\ & = & \cyclic\,\frac{1}{16}\bigl(\epsilon_k\langle R^\nabla_{X,Y}(\nab{Z}{J})e_k,\overline{e_k}\rangle+ 
\epsilon_k\langle R^\nabla_{X,Y}(\nab{Z}{J})\overline{e_k},e_k\rangle\bigr)
\end{eqnarray*}
We are going to compute $\dx\omega(u,v,\overline{z})$ for any $u,v,z\in\XIS^+$, the $+i$-eigenspace of $J$, because it corresponds to the computation of $\dx\omega^{2,1}$ (or $\dx\omega^{1,2}$ by conjugation of the real form).

The integrability condition implies $(\nab{u}{J})v=0,\ \forall u,v\in\XIS^+$ because $\nab{u}{v}\in\XIS^+$. Of course, we have $(\nab{\overline{u}}{J})\overline{v}=0$ too. Let $\xi$ denote any index and let $\nab{\xi}{e_k}=\gamma_{\xi,k}^he_h+\gamma_{\xi,k}^{\overline{h}}\overline{e_h}$. Then
\[  (\nab{\xi}{J})e_k=(i-J)\nab{\xi}{e_k}=2i\gamma_{\xi,k}^{\overline{h}}\overline{e_h} \qquad(\mbox{hence}\quad  \gamma_{j,k}^{\overline{h}}=0 ) \]
and
\[  (\nab{\xi}{J})\overline{e_k}=(-i-J)\nab{\xi}{\overline{e_k}}=-2i\gamma_{\xi,\overline{k}}^he_h
\qquad(\mbox{hence}\quad  \gamma_{\overline{\xi},k}^{\overline{h}}=\overline{\gamma_{\xi,\overline{k}}^h} ). \]
Notice $\nabla J$ permmutes the $+$ and $-$ $i$-eigenspaces. From
\[   \langle \gamma_{\xi,k}^{\overline{h}}\overline{e_h},e_j\rangle= \langle\nab{\xi}{e_k},e_j\rangle=-\langle e_k,\nab{\xi}{e_j}\rangle=-\langle e_k, \gamma_{\xi,j}^{\overline{h}}\overline{e_h}\rangle  \]
we find $\epsilon_j\gamma_{\xi,k}^{\overline{j}}=-\epsilon_k\gamma_{\xi,j}^{\overline{k}}$. Finally,
\[ \dx\omega(u,v,\overline{z}) = \frac{1}{16}\bigl(\epsilon_k\langle R^\nabla_{u,v}(\nab{\overline{z}}{J})e_k,\overline{e_k}\rangle+ 
\epsilon_k\langle R^\nabla_{\overline{z},u}(\nab{v}{J})\overline{e_k},e_k\rangle+ 
\epsilon_k\langle R^\nabla_{v,\overline{z}}(\nab{u}{J})\overline{e_k},e_k\rangle\bigr) . \]
But using the symmetries $\langle R_{u,v}a,b\rangle=\langle R_{a,b}u,v\rangle=-\langle R_{u,v}b,a\rangle$, we find
\[  \langle R^\nabla_{v,\overline{z}}(\nab{u}{J})\overline{e_k},e_k\rangle= -2i\gamma_{u,\overline{k}}^h\langle R^\nabla_{v,\overline{z}}e_h,e_k\rangle=
-2i\gamma_{u,\overline{k}}^h\langle R^\nabla_{e_h,e_k}v,\overline{z}\rangle=0  \]
and therefore we may continue from above
\begin{eqnarray*}
\dx\omega(u,v,\overline{z}) & = & \frac{1}{16}\epsilon_k\langle R^\nabla_{u,v}(\nab{\overline{z}}{J})e_k,\overline{e_k}\rangle  \\
& = & \frac{i}{8}\gamma_{\overline{z},k}^{\overline{h}}
\epsilon_k\langle R^\nabla_{u,v}\overline{e_h},\overline{e_k}\rangle  \\
& = & \frac{i}{8}\gamma_{\overline{z},k}^{\overline{h}} \epsilon_k \bigl(
\langle v,\overline{e_h}\rangle\langle u,\overline{e_k}\rangle-\langle u,\overline{e_h}\rangle\langle v,\overline{e_k}\rangle \bigr) .
\end{eqnarray*}
Now we apply this to $u=e_\alpha,\ v=e_\beta$. We get
\[   \dx\omega(e_\alpha,e_\beta,\overline{z})= \frac{i}{2}(\epsilon_\beta\gamma_{\overline{z},\alpha}^{\overline{\beta}}-\epsilon_\alpha\gamma_{\overline{z},\beta}^{\overline{\alpha}}) =i\epsilon_\beta\gamma_{\overline{z},\alpha}^{\overline{\beta}} .\]
On the other hand, using formula (\ref{formula simples e estranha}) we immediately find
\begin{eqnarray*}
\dx\omega(e_\alpha,e_\beta,\overline{z}) & = & \langle(\nab{e_\alpha}{J})e_\beta,\overline{z}\rangle + \langle(\nab{e_\beta}{J})\overline{z},e_\alpha\rangle + \langle(\nab{\overline{z}}{J})e_\alpha,e_\beta\rangle \\
& = & \langle(\nab{\overline{z}}{J})e_\alpha,e_\beta\rangle \ =\ 2i\gamma_{\overline{z},\alpha}^{\overline{h}}\langle\overline{e_h},e_\beta\rangle \\ 
& = & 4i\epsilon_\beta\gamma_{\overline{z},\alpha}^{\overline{\beta}}.
\end{eqnarray*}
This implies $\dx\omega=0$.
\end{proof}
\begin{teo}
There is no integrable orthogonal complex structure on $S^{2n}_{2q}$.
\end{teo}
\begin{proof}
Such a complex structure would have to be pseudo-K\"ahlerian in contradiction with proposition \ref{naoe'psdKahler}.
\end{proof}
\ \\
\textbf{Remark.}
$S^6_4$ does not admit an orthogonal integrable complex structure, but it has a nearly pseudo-K\"ahler structure with respect to the usual metric. In fact we can generalize E.~Calabi's construction as follows. We first consider $\R^3$ with a Lorentz metric $g$ and let $(e_1,e_2,e_3)$ denote an orthonormal basis with signature $+--$. Then a cross pruduct is well defined by $g(u\times v,w)= \vol (u,v,w)$, where the $\vol = e^{(123)}$, --- which can be extended to elements of $R^4$; writing $a=(a_0,a'),\ b=(b_0, b')$, then $a\times b=-a_0b' +b_0a' + a' \times b'$ and a quaternionic multiplication can be given as 
\[ a\cdot b= (a_0b_0-g(a',b'), a_0b' +b_0a' + a' \times b' ) . \]
Thus 
\[ a\times b= \pig (\overline{b}\cdot a)   \]
where $\overline{b}= (b_0, -b')$ is the conjugate. Then $\R^4$ adopts the signature $++--$ and we can define a new fixed metric on $\R^8= \R^4 \times \R^4$ with signature $++--\,++--$.

The definition of a cross product as in the Cayley-Dickson process is then possible: letting $u=(a, \alpha),\  v=(b, \beta)\in\R^8$, 
\[ u \times v = (a\times b - \alpha\times \beta , (\alpha\cdot\overline{\beta}- \beta\cdot\overline{\alpha})) .  \]
Now we take the pseudo-sphere $S=S^6_4=\{x\in \R^7: g(x,x)=1\} \subset  0 \times \R^7 \subset \R^4 \times \R^4$.
Since $\R^7$ has signature $+-- ++--$, this implies $S$ with signature $-- ++--$. Finally, if $x\in S$ and $u\in T_xS$, then the map defined by $J_x(u)= x \times u$ is an orthogonal almost complex $J$. One proves this $J$ is nearly pseudo-K\"ahler and non-integrable, just as in the Riemannian case.

\newpage

\bibliographystyle{plain}

\begin{thebibliography}{10}


\bibitem{Albu}
R.~Albuquerque.
\newblock {\em On the Twistor Theory of Symplectic Manifolds}.
\newblock PhD thesis, Mathematics Institute, University of Warwick, 2002.

\bibitem{AR}
R.~Albuquerque and J.~Rawnsley.
\newblock {\em Twistor Theory of Symplectic Manifolds}.
\newblock to appear in Jour. of Geom. and Physics.

\bibitem{Atiyah}
M.~Atiyah, N.~Hitchin and I.~Singer.
\newblock {\em Self-duality in four-dimensional riemannian geometry}.
\newblock Proc. Royal Soc. of London, A 362(1711):425--461, 1978.

\bibitem{BR}
M.~Barros and A.~Romero.
\newblock {\em Indefinite K\"{a}hler Manifolds}.
\newblock Math. Annalen, 261:55--62, 1982.

\bibitem{Besse}
A.L.~Besse.
\newblock {\em Einstein Manifolds}.
\newblock Springer, Berlin, 1987.

\bibitem{Blair}
D.~E.~Blair, J.~Davidov, O.~Muskarov.
\newblock {\em Isotropic K\"ahler hyperbolic twistor spaces}.
\newblock Journal of Geometry and Physics, 52, pages 74--88, 2004.

\bibitem{Bu}
J.~B.~Butruille\ss.
\newblock {\em Espace de twisteurs r\'eduit d'une vari\'et\'e presque complexe de dimension 6}.
\newblock DG/0503150.

\bibitem{Hel}
S.~Helgason.
\newblock {\em Differential Geometry, Lie Groups, and Symmetric Spaces}.
\newblock Academic Press, 1978.

\bibitem{Hit}
N.~J. Hitchin.
\newblock {\em K\"ahlerian twistor spaces}.
\newblock Proc. London Math. Soc. (3), 43(1):133--150, 1981.

\bibitem{KobNomi}
S.~Kobayashi and K.~Nomizu.
\newblock {\em Foundations of Differential Geometry}.
\newblock Vol. 1 and 2, Wiley Classics Library, edition 1996.

\bibitem{Lebr}
C.~Lebrun.
\newblock {\em Orthogonal Complex Structures on $S^6$}.
\newblock Proc. of the Amer. Math. Soc., Vol 101, 1, September 1987.

\bibitem{Obri}
N.~O'Brian and J.~Rawnsley.
\newblock {\em Twistor spaces}.
\newblock Ann. of Global Analysis and Geometry, 3(1):29--58, 1985.

\bibitem{Rawnsley}
J.~Rawnsley.
\newblock {\em $f$-structures, $f$-twistor spaces and harmonic maps}.
\newblock Springer, Berlin, pages 85--159, 1985.
\newblock Geometry seminar ``Luigi Bianchi'' II.

\bibitem{Sal}
S.~Salamon.
\newblock {\em Harmonic and holomorphic maps}.
\newblock Springer, Berlin, pages 161--224, 1985.
\newblock Geometry seminar ``Luigi Bianchi'' II.

\bibitem{Vais}
I.~Vaisman.
\newblock {\em Symplectic twistor spaces}.
\newblock Jour. of Geom. and Phys., 3(4):507--524, 1986.

\end{thebibliography}

\end{document}